\documentclass[12pt,leqno]{article}
\usepackage[francais]{babel}
\usepackage[latin1]{inputenc}
\usepackage[T1]{fontenc}
\usepackage{amsmath,amsfonts,amssymb,amsthm,amscd}
\usepackage{bm}
\date{}

\newcommand{\ba}{\backslash}
\newcommand\rg{\rightarrow}

\newcommand\Z{\mathbb{Z}}

\newcommand\C{\mathbb{C}}

\newcommand\N{\mathbb{N}}

\newcommand\R{\mathbb{R}}

\newcommand\Q{\mathbb{Q}}

\newcommand\Gal{\mathrm{Gal}}

\newcommand{\sumn}{\sum\limits_{n=1}^{+\infty}}
\newcommand{\sumk}{\sum\limits_{k=0}^{+\infty}}

\newcommand{\newsmile}{\smallsmile_{\hspace{-3mm}\displaystyle\smallfrown}}

\def\adots{\mathinner{\mkern1mu\raise1pt\vbox{\kern7pt\hbox{.}}
\mkern2mu\raise4pt\hbox{.}
\mkern2mu\raise7pt\hbox{.}\mkern1mu}}

\begin{document}

\title{Fonctions complètement multiplicatives\\de somme nulle}

\author{Jean--Pierre Kahane et Eric Saias}

\date{juillet 2015}

\maketitle

\noindent\textbf{Abstract. Completely multiplicative functions whose sum is zero ($CMO$)}
\vspace{2mm}

The paper deals with $CMO$, meaning completely multiplicative ($CM$) functions $f$ such that $f(1)=1$ and $\sum\limits_1^\infty f(n)=0$. $CM$ means $f(ab)=f(a)f(b)$ for all $(a,b)\in \N^{*2}$, therefore $f$ is well defined by the $f(p)$, $p$ prime. Assuming that $f$ is $CM$, give conditions on the $f(p)$, either necessary or sufficient, both is possible, for $f$ being $CMO$ : that is the general purpose of the authors.

The $CMO$ character of $f$ is invariant under slight modifications of the sequence $(f(p))$ (theorem~3). The same idea applies also in a more general context (theorem~4).

After general statements of that sort, including examples of $CMO$ (theorem~5), the paper is devoted  to ``small'' functions, that is, functions of the form $\frac{f(n)}{n}$, where the $f(n)$ are bounded. Here is a typical result : if $|f(p)|\le 1$ and $Re\, f(p)\le0$ for all $p$, a necessary and sufficient condition for $\big(\frac{f(n)}{n}\big)$ to be $CMO$ is $\sum \, Re\, f(p)/p=-\infty$ (theorem~8). Another necessary and sufficient condition is given under the assumption that $|1+f(p)|\le 1$ and $f(2)\not=-2$ (theorem~7). A third result gives only a sufficient condition (theorem~9). The three results apply to the particular case $f(p)=-1$, the historical example of Euler.

Theorems 7 and 8 need auxiliary results, coming either from the existing literature (Hal\'asz, Montgomery--Vaughan), or from improved versions of classical results (Ingham, Ska\l ba) about $f(n)$ under assumptions on the $f*1(n)$, * denoting the multiplicative convolution (theorems~10~and~11).

\vspace{4mm}


\section{Les fonctions $CMO$, théorème 1}

Appelons  fonction complètement multiplicative de somme nulle (ou fonction $CMO$) toute fonction $f$ de $\N^*$ dans $\C$ non identiquement nulle, et qui vérifie les deux conditions suivantes :
$$
f(ab) = f(a)f(b) \qquad \mathrm{pour\ tout\ } (a,b)\in \N^{*2}
$$
et
$$
\sum_{n=1}^{+\infty}f(n)=0\,.
$$
Précisons que conformément à l'usage, cette dernière  identité signifie que la suite des sommes partielles $\Big(\sum\limits_{n=1}^N f(n)\Big)_{N\ge 1}$ converge dans $\C$ quand $N$ tend vers $+\infty$, et que sa limite est nulle.

\vspace{2mm}
L'objet de ce présent travail est de donner les informations que nous avons pu glaner sur l'ensemble $E$ des fonctions $CMO$. Déjà, on a~le

\vspace{4mm}

\noindent \textbf{Théorème 1}.--- \textit{L'ensemble $E$ est non vide.}

\vspace{4mm}

Nous engageons le lecteur à chercher à démontrer ce résultat avant de passer à la suite.
\vspace{2mm}

Les remarques de base que sont les théorèmes 2 et 3  seront démontrés au fil de la plume au paragraphe~2. La preuve du théorème~4, qui est dans le même esprit que celle du théorème~3, sera laissée au soin du lecteur.  Le paragraphe~3 est relatif à la notion de série de Dirichlet génératrice. Le théorème~5, qui permet de valider au passage le théorème~1, sera établi aux paragraphes~4 et 5. Le paragraphe~6 est dévolu principalement à quelques remarques générales sur les fonctions $CMO$. Le théorème~6  sera démontré au paragraphe~7. Les théorèmes 10 et 11 sont énoncés au paragraphe~8, et seront démontrés respectivement aux paragraphes~10 et 11. Les théorèmes~7, 8 et 9 seront démontrés respectivement aux paragraphes~12, 13 et 14. Enfin au paragraphe~9, on évoque un lien  avec   les zéros des fonctions $L$ d'Artin.

Signalons que dans tout ce texte, on utilisera la lettre $n$ pour désigner un entier générique $\ge 1$, et la lettre $p$ pour désigner un nombre premier générique.

\section{Les phénomènes généraux de base, théorèmes 2, 3 et 4}

On a envie d'emblée de voir un exemple concret de fonction $CMO$. Comme une fonction complètement multiplicative (ou fonction $CM$) est déterminée par ses valeurs aux nombres premiers, la première idée pour construire un tel exemple consiste à choisir une suite $f(p)$ sur les nombres premiers, et d'ajuster ces valeurs de telle sorte que $\sum\limits_{n=1}^{+\infty}f(n)=0$. Le plus simple pour assurer la convergence de la série est de demander sa convergence absolue. Mézalor la formule du produit Eulérien
$$
\sum_n f(n) =\prod _p \frac{1}{1-f(p)}
$$
empêche la série de s'annuler !
\begin{equation}
\mathrm{Une\ fonction\ } CMO\ \textrm{\ n'est\ pas\ sommable}\,.
\end{equation}

\noindent Une fonction $CMO$ n'est pas sommable. Ce fait appelle plusieurs remarques, dont les trois premières montrent son caractère optimal dans des directions différentes.

\vspace{2mm} 
\noindent\textbf{1)} Il est facile de montrer que si $S\in \C^*$, il existe une fonction $CM$ $f$ telle que
$$
\sum_{n\ge 1}|f(n)| <+\infty\qquad \mathrm{et}\qquad \sum_{n\ge 1} f(n)=S\,.
$$
Autrement dit, zéro est la seule valeur interdite pour la somme d'une fonction $CM$ et sommable.

\vspace{2mm} 
\noindent\textbf{2)} Il est facile de déduire du théorème~8 du présent travail que l'on peut trouver des fonctions $f$ $CMO$ avec la fonction positive de $x$, $\sum\limits_{2\le n\le x} |f(n)|$, qui croît arbitrairement lentement vers $+\infty$.

\vspace{2mm} 
\noindent\textbf{3)} L'assertion (2.1) est mise en défaut si on suppose uniquement $f$ multiplicative. Par exemple, si $f(1)=-f(2)=1$ et $f$ est nulle partout ailleurs, alors $f$ est multiplicative, sommable et de somme nulle.

\vspace{2mm} 
\noindent\textbf{4)} Il résulte en particulier de (2.1) que pour tout $s\in \C$ avec $Re\ s >1$, $1/n^s$ n'est pas $CMO$.
C'est l'argumentation classique dûe à Riemann \cite{Rie} pour montrer que $\zeta(s)$ ne s'annule pas dans le demi--plan $Re\ s>1$.

\vspace{2mm} 
\noindent\textbf{5)} Plus généralement, il est à noter qu'aucune fonction $1/n^s$ n'est $CMO$, car pour $Re\ s \le 1$ la série $\sum\limits_{n=1}^{+\infty} 1/n^s$ diverge. Or on sait que les fonctions $1/n^s$ sont les seules fonctions complètement multiplicatives à être régulières en un certain sens. Par exemple il résulte du théorème~2 de l'élégant travail de Wirsing et Zagier \cite{WiZa}, que ce sont les seules pour lesquelles $f(n+1)/f(n)$ tend vers $1$. En raison donc de l'irrégularité de la suite $(f(n))_{n\ge 1}$ quand $n$ parcourt en croissant l'ensemble des entiers, irrégularité qui est liée à la dépendance de $f(n)$ à la décomposition de $n$ en facteurs premiers, on peut dire que les fonctions $CMO$ sont de nature arithmétique.

\vspace{2mm} 
\noindent\textbf{6)} D'après P\'olya et Szegö (voir le problème 129 du \S~1 de \cite{PoSz}), l'assertion~(2.1) remonte au moins à Haar, sous la forme plus générale suivante : soit $f$ une fonction arithmétique non nulle vérifiant $\sum\limits_{n=1}^{+\infty} f(kn)=0$ pour tout entier $k\ge1$, alors $f$ n'est pas sommable.

\vspace{2mm} 
\noindent\textbf{7)} Soit $f$ une fonction $CMO$. On a en particulier que $(f(n))_{n\ge 2}$ et ces sous--suites $\big(f(n^k)\big)_{k\ge1}=\big((f(n))^k\big)_{k\ge 1}$ tendent vers~$0$. On en déduit~que
$$
\forall n\ge 2, \quad |f(n)| <1\,.
$$
De plus, le supremum des $|f(n)|$, comme celui des $|f(p)|$, est atteint. Enfin, en réutilisant le caractère $CM$, on observe que ces deux nombres sont égaux. On a donc montré~que
\begin{equation}
\sup_{n\ge 2}|f(n)|=\max_{n\ge 2}|f(n)| =\max_p |f(p)| = \sup_p |f(p)|<1\,.
\end{equation}

\vspace{2mm} 

\noindent\textbf{8)} On a vu que si $f$ est $CMO$, $\sum\limits_n|f(n)|=+\infty$. En fait, on a même  $\sum\limits_p|f(p)|=+\infty$. En effet si  $\sum\limits_p |f(p)|<+\infty$, en utilisant (2.2) et le produit Eulérien on aurait
$$
\sum_n f(n) =\prod_p (1-f(p))^{-1} \not=0\,.
$$
Récapitulons. Nous avons montré le

\vspace{2mm}

\noindent \textbf{Théorème 2}.--- \textit{Soit $f$ une fonction $CMO$. Alors}
\begin{equation}
\sup_{n\ge 2}|f(n)| = \max_{n\ge 2}|f(n)| = \max_p |f(p)|=\sup_p |f(p)|<1
\end{equation}
\textit{et}
\begin{equation}
\sum_p|f(p)| = +\infty\,.
\end{equation}

On s'intéresse maintenant à la stabilité du caractère $CMO$. Dans quelle mesure une fonction $CM$, proche d'une fonction $CMO$, est encore $CMO$ ? Le résultat suivant donne une réponse générale.

\vspace{2mm}

\noindent \textbf{Théorème 3}.--- 

\vspace{2mm}

\begin{tabular}{rcll}
\textit{Soient} & $f$ &\textit{une fonction} &$CM$\\
\textit{et}	& $g$ &\textit{une fonction} & $CMO$
\end{tabular}

\noindent\textit{vérifiant}
\begin{equation}
\forall p,\quad |f(p)| < 1
\end{equation}
\textit{et}
\begin{equation}
\sum_p |(f-g)(p)| <+\infty\,.
\end{equation}

\textit{Alors la fonction $f$ est $CMO$.}

\vspace{2mm}

\noindent\textbf{Démonstration}. En effet notons $F(s) = \sumn f(n)/n^s$ et $G(s)=\sumn g(n)/n^s$. Alors les hypothèses~(2.5) et (2.6) permettent de montrer que la\break série de Dirichlet $(F/G)(s)$ converge absolument en $s=0$. Donc\break $F(0)=(F/G)(0)\cdot G(0)=0$, d'après le caractère $CMO$ de $g$. Autrement dit, $f$ est~$CMO$.

Nous serons amenés à travailler dans le cadre plus général où l'une des deux fonctions est multiplicative, sans être nécessairement complètement multiplicative. Nous laissons le soin au lecteur de vérifier que la même démarche que  celle de la preuve du théorème~3, permet d'établir le résultat suivant

\vspace{4mm}

\noindent \textbf{Théorème 4}.--- \textit{Soient $f$ une fonction complètement multiplicative, et $g$ une fonction multiplicative vérifiant}
\begin{eqnarray}
&\forall p,\ \ |f(p)|<1\ \hbox to 3,5cm\\
&\forall p,\ \sumk g(p^k) \in \C^*\hbox to 3,5cm\\\
&\sum\limits_p |(1-f(p)) \sumk g(p^k) -1| <+\infty\,.		 
\end{eqnarray}
\textit{Alors les deux assertions (2.10) et (2.11) suivantes sont équivalentes}
\begin{eqnarray}
\sumn f(n)=0\\
\sumn g(n)=0\,.
\end{eqnarray}

\vspace{2mm}

\section{Série de Dirichlet génératrice associée à une fonction $CM$}

Soit $f$ une fonction complètement multiplicative. On remarque que $f$ est $CMO$ si et seulement si sa série de Dirichlet génératrice associée
\setcounter{equation}{0}
\begin{equation}
F(s)= \sumn \frac{f(n)}{n^s}
\end{equation}
converge en $s=0$, et y est de somme nulle. En réalité,  tous les résultats de ce présent travail résultent  de l'étude des séries de Dirichlet à coefficients complètement multiplicatifs. Dorénavant,  quand $f(n)$ désignera une fonction $CM$, on utilisera systématiquement la notation (3.1) pour sa série de Dirichlet associée. On notera également dans tout ce travail
$$
s=\sigma +it\,.
$$

\section{Les exemples de Dirichlet--Riemann de\\ fonction $CMO$}

Soit $f$ une fonction $CM$. Notons $\sigma_c$ l'abscisse de convergence de la série de Dirichlet associée $F(s)$. Le cas le plus simple d'existence de fonction $CMO$ est quand la fonction $F(s)$ s'annule dans le demi--plan $\{Re\ s >\sigma_c\}$. C'est le cas pour
\setcounter{equation}{0}
\begin{equation}
L_\chi(s)= \sumn \frac{\chi(n)}{n^s}
\end{equation}
quand $\chi$ est un caractère de Dirichlet non principal. En effet dans ce cas on sait que $\sigma_c=0$ et que la fonction $L_\chi(s)$ s'annule une infinité de fois dans le demi--plan $\{Re\ s >0\}$. Pour tout tel zéro $\rho$, on a donc $\chi(n) /n^\rho$ qui est~$CMO$.

D'un point de vue historique, ces exemples de fonctions $CMO$ ont pu voir le jour grâce aux travaux successifs de Dirichlet et de Riemann. C'est en effet Dirichlet qui dans son travail magnifique \cite{Dir} de 1837 introduit les caractères $\chi$ modulo $q$ et leur fonction $L_\chi(s)$ associée. Mais il ne considérait ces fonctions $L_\chi(s)$ que pour la variable $s$ réelle. C'est Riemann dans son célèbre article \cite{Rie} de 1859, qui a l'idée de considérer la fonction $\zeta(s)$ d'Euler comme une fonction de variable complexe. De plus il démontre que la fonction $\zeta(s)$ admet une infinité de zéros dans la bande critique $\{0\le Re \ s\le 1\}$. Avec l'équation fonctionnelle, on en déduit plus précisément que la fonction $\zeta(s)$ s'annule une infinité de fois dans la demi--bande $\{\frac{1}{2}\le Re\ s \le 1\}$.

C'est en travaillant de manière analogue avec $L_\chi(s)$ où $\chi$ est un caractère de Dirichlet non principal, que l'on obtient que la formule (4.1) définit une fonction sur le demi--plan $Re\ s >0$ du plan complexe, qui s'annule une infinité de fois dans la demi--bande $\{\frac{1}{2}\le Re\ s \le 1\}$.

On appellera dorénavant exemple de fonction $CMO$ de Dirichlet--Riemann toute fonction $f(n)=\chi(n)/n^\rho$ où $\chi$ est un caractère de Dirichlet non principal et $\rho$ un nombre complexe vérifiant $Re\ \rho>0$ et $L_\chi(\rho)=0$.

\section{L'exemple d'Euler}
\setcounter{equation}{0}
Considérons à présent la fonction de Liouville $\lambda(n)$ qui est la fonction $CM$ définie sur les nombres premiers par $\lambda(p)=-1$. Ainsi $\lambda(n)=\pm1$.

Commençons par calculer sa série de Dirichlet associée. Pour $Re\ s>1$, la série $\sumn \lambda(n)/n^s$ est absolument convergente et~on~a
\begin{equation}
\sum_n \frac{\lambda(n)}{n^s} = \prod_p \frac{1}{1-\frac{\lambda(p)}{p^s}} = \prod_p \frac{1}{1+\frac{1}{p^s}}=\prod_p \frac{1-\frac{1}{p^s}}{1-\frac{1}{p^{2s}}}=\frac{\zeta(2s)}{\zeta(s)}\,.
\end{equation}
Or la fonction $\zeta(2s)/\zeta(s)$ s'annule en $s=1$. Si on savait que l'abscisse de convergence de $\sumn \lambda(n)/n^s$ est strictement inférieure à $1$, on pourrait conclure immédiatement que $\lambda(n)/n$ est $CMO$. Mais on est très loin de le savoir actuellement. Notons de plus qu'il résulte de (5.1) que l'on~a
$$
\lim_{\sigma\searrow1} \sumn \frac{\lambda(n)}{n^\sigma}= \lim_{\sigma\searrow1} \frac{\zeta(2\sigma)}{\zeta(\sigma)}=0\,,
$$
ce qui exprime que la série $\sum\limits_n \lambda(n)/n$ a pour somme $0$ quand on lui applique un procédé de sommation convenable. Mais c'est bien insuffisant pour montrer la convergence de la série au sens de la convergence des sommes partielles. Indiquons à présent une preuve de l'identité $\sumn \lambda(n)/n=0$.

\vspace{2mm}

\textbf{Preuve de $\sumn \lambda(n)/n=0$ par l'analyse de Fourier}

\vspace{2mm}
Soit $x\ge 2\ge 2a>0$. On part de la formule
\begin{equation}
\int_\R \frac{\sin xt}{t} \Big(\frac{\sin at}{at}\Big)^2\ e^{-itu}dt := \pi \, a\, w_{x,a}(u)
\end{equation}
où
$$
w_{x,a}(u) \left |
\begin{array}{lll}
=1 &\mathrm{pour} &|u|\le x-2a\\
=0 &\mathrm{pour} &|u|\ge x+2a\\
\in [0,1] &\mathrm{pour} &u\in \R\,.
\end{array}
\right.
$$
Ainsi, pour $\sigma>1$,
$$
\int_\R \frac{\zeta(2\sigma+2it)}{\zeta(\sigma+it)} \frac{\sin xt}{t} \Big(\frac{\sin at}{at}\Big)^2 dt = \pi\ a \sum_n \frac{\lambda(n)}{n^\sigma} w_{x,a}(\log n)\,.
$$
Or la fonction
$$
\frac{\zeta(2\sigma+2it)}{\zeta(\sigma+it)} \frac{\sin xt}{t}\Big(\frac{\sin at}{at}\Big)^2
$$
tend vers la fonction
$$
\frac{\zeta(2+2it)}{\zeta(1+it)} \frac{\sin xt}{t}\Big(\frac{\sin at}{at}\Big)^2
$$
dans $L^1(\R)$ quand $\sigma\searrow 1$, d'après les propriétés de la fonction $\zeta(s)$ ci--dessous (cf. théorèmes~II.3.1 et II.3.22 de \cite{Ten})
\begin{eqnarray}
\zeta^{-1}(s) \newsmile s-1\hspace{1,5cm} \qquad \qquad (|s-1|\le 2)\\
\zeta^{-1}(s) \ll \log\, 3|t| \qquad\qquad(\sigma\ge 1,\ |t|\ge 1/2)\,.
\end{eqnarray}
Il en résulte que
$$
\sum_{\log n\le x} \frac{\lambda(n)}{n} = \frac{1}{\pi a}\int_\R \frac{\zeta(2+2it)}{\zeta(1+it)}\Big(\frac{\sin at}{at}\Big)^2 \sin\, xt\ dt +O(a+e^{-x})
$$
où la constante implicite dans le symbole $O$ est absolue. On conclut grâce au lemme de Riemann--Lebesgue (a fixé) puis en faisant tendre $a$ vers~$0$.

Il est à noter que la formule $\sum\limits_{n=1}^{+\infty} \lambda(n)/n=0$ que nous venons de démontrer, se trouve en un sens faible dans le Theorema~18 du travail fondamental \cite{Eul} d'Euler en 1737.

Signalons que l'on peut prouver de la même manière que
\begin{equation}
\sumn \frac{\mu(n)}{n}=0
\end{equation}
où $\mu(n)$ désigne la fonction de Möbius. Il suffit de remplacer dans notre argumentation, la formule (5.1) par la formule légèrement plus simple :
$$
\sum_n \frac{\mu(n)}{n^s} = \frac{1}{\zeta(s)}\qquad (Re\ s >1)\,.
$$
On peut alors déduire le caractère $CMO$ de $\lambda(n)/n$ de la formule~(5.5), en appliquant le théorème~4.

\section{Récapitulation, théorème 5}

Récapitulons dans un énoncé formel les exemples de fonction $CMO$ que nous avons mis en lumière dans les deux paragraphes précédents.

\vspace{2mm}

\noindent \textbf{Théorème 5} (exemples d'Euler et de Dirichlet--Riemann de fonctions $CMO$)--- 

\vspace{-1mm}
\textit{(i) La fonction $\lambda(n)/n$ est $CMO$.}

\textit{(ii) Soient $\chi$ un caractère de Dirichlet non principal et $\rho$ un nombre complexe tels que $Re\ \rho >0$ et $L_\chi(\rho)=0$. Alors la fonction $\chi(n)/n^\rho$ est~$CMO$.
}
\vspace{2mm}

Nous avons réuni l'exemple d'Euler et les exemples de Dirichlet--Riemann en un seul énoncé car ce sont des exemples de fonctions $CMO$ déjà connus. Il est à noter cependant qu'ils sont de nature très différente.

Pour $\lambda(n)/n$, la série de Dirichlet associée $F(s)= \sumn \lambda(n)/n^s=\zeta (2s)/\zeta(s)$ ne s'annule qu'en un seul point, bien spécifié, $s=1$. Ce point est situé sur la droite d'abscisse de convergence absolue de $F(s)$ ; cela permet de disposer du produit Eulérien de $F(s)$ dans tout le demi--plan $Re\ s>1$, et plus particulièrement au voisinage à droite de la droite~$Re\,s=1$.

Les exemples de Dirichlet--Riemann proviennent de l'étude des séries $L_\chi(s)$ de Dirichlet associées à des caractères de Dirichlet $\chi$ non principaux. C'est l'équation fonctionnelle de $L_\chi(s)$ et le principe de l'argument qui permet de montrer que $L_\chi(s)$ admet une infinité de zéros  dans la bande critique fermée $0\le Re\ s\le 1$. En fait, on sait qu'ils sont tous dans la bande critique ouverte $0<Re\ s<1$, ce qui permet à chacun de ces zéros de fournir un exemple de fonction $CMO$. Contrairement à l'exemple d'Euler, tous ces zéros sont strictement à gauche de la droite d'abscisse de convergence absolue de $L_\chi(s)$. Enfin, il est à noter que ces zéros $\rho$ (et donc aussi les fonctions $CMO$ $\chi(n)/n^\rho$ qui leur correspondent) ne sont pas connus de manière précise dans la mesure où on ne connait pas de formule qui donne la valeur exacte de leur partie imaginaire.

\section{Petites fonctions $CMO$, théorèmes 6, 7, 8 et~9}
\setcounter{equation}{0}

Etudier les fonctions $CMO$ générales nous semble être un problème difficile. Dans le cas particulier de l'exemple d'Euler $\lambda(n)/n$, on a utilisé son produit Eulérien
$$
\sum_n \frac{\lambda(n)}{n^s} = \prod_p \frac{1}{1+1/p^s}
$$
qui est valable dans tout le demi--plan vertical $\{Re\, s>1\}$,  au voisinage à droite du point $s=1$ qui nous intéresse.

Appelons ici petite fonction $CM$ (respectivement $CMO$) toute fonction de la forme $f(n)/n$ où $f$ est une fonction $CM$ ($\mathrm{resp}. CMO$) bornée sur les nombres premiers. Le point  est que l'on a encore dans ce cadre l'outil fondamental qu'est son produit Eulérien
\begin{equation}
F(s) = \sum_n \frac{f(n)}{n^s}=\prod_p \frac{1}{1-f(p)/p^s}\,,\quad (Re\ s>1)\,.
\end{equation}

Si $f(n)/n$ est $CMO$, on sait grâce au théorème~2 que $\sum\limits_p |f(p)|/p=+\infty$. Commençons par remarquer qu'en utilisant le produit d'Euler (7.1), on peut dire un peu plus pour les petites fonctions $CMO$. En effet sur la demi--droite $\sigma>1$, on~a
\begin{equation}
\!\!\!\!\Big|\!\!\sum_n \!\frac{f(n)}{n^\sigma}\Big| \!\!=\!\!\Big|\!\prod_p \frac{1}{1\!\!-\!\!f(p)/p^\sigma}\Big| \!\!\newsmile \!\!
\exp \Big\{\!\sum_p \!\!\frac{Re(f(p))}{p^\sigma}\Big\} \!\!\ge\!\!\exp\!\Big\{\!\!\!\!\!\sum_{Re\,f(p)<0}\!\!\!\!\!\!\frac{Re\,f(p)}{p^\sigma}\Big\}.
\end{equation}
Or en transposant aux séries de Dirichlet le théorème d'Abel sur les séries entières, on~a
$$
\lim_{\sigma\searrow 1} \sum_n \frac{f(n)}{n^\sigma} = \sumn \frac{f(n)}{n}=0\,.
$$
Donc en passant à la limite dans l'inégalité (7.2), on obtient par convergence monotone
$$
0\ge \exp \Big\{\sum_{Re\, f(p)<0} \frac{Re\,f(p)}{p}\Big\}\,,
$$
d'où finalement
$$
\sum_{Re\,f(p)<0} \frac{Re\,f(p)}{p} = -\infty\,.
$$
On a donc établi le

\vspace{2mm}

\noindent \textbf{Théorème 6}.--- \textit{Soit $f(n)/n$ une petite fonction $CMO$. On a alors}
$$
\sum_{Re\,f(p)<0} \frac{Re\,f(p)}{p}= -\infty\,.
$$

\vspace{1mm}

Revenons maintenant à l'exemple fondamental d'Euler de petite fonction $CMO$ : $\lambda(n)/n$. On le sait, il est facile de déduire les deux formules suivantes l'une de l'autre (cf.~théorème~4)
$$
\sumn \frac{\lambda(n)}{n}=0
$$
et 
\begin{equation}
\sumn \frac{\mu(n)}{n} = 0\,.
\end{equation}
A l'aide d'une inversion de Möbius qui utilise le principe de l'hyperbole de Dirichlet, on peut généraliser cette dernière identité de la manière suivante (théorème~II.7.24 de~\cite{Ten}). Dans tout ce présent travail, le symbole $*$ désigne l'opération de convolution multiplicative.

\vspace{2mm}

\noindent \textbf{Théorème $I+S$} (Ingham + Ska\l ba).--- \textit{Soient $f$ et $g$ deux fonctions arithmétiques, liées par la relation}
$$
g=f * 1\,.
$$
\textit{On suppose que $g$ est bornée et que}
\begin{equation}
\lim_{x\rg +\infty} \frac{1}{x} \sum_{n\le x} g(n) =0\,.
\end{equation}
\textit{On a alors}
$$
\sumn \frac{f(n)}{n}=0\,.
$$

\vspace{2mm}

Le cas particulier (7.3) correspond à $g=\delta_1$.

Or grâce à l'important travail de Wirsing de 1967 \cite{Wir} et sa généralisation l'année suivante par Hal\'asz (\cite{Hal} ou \cite{Ten}), on sait caractériser les fonctions multiplicatives $g$ à valeurs dans le disque unité qui vérifient (7.4). Cela permet de donner une condition suffisante pour qu'une petite fonction $CM$, $f(n)/n$, proche en un certain sens de l'exemple
paradigmatique d'Euler $\lambda(n)/n$, soit~$CMO$.

En travaillant un peu plus, nous avons obtenu le résultat suivant :

\vspace{2mm}

\noindent \textbf{Théorème 7}.--- \textit{Soit $f$ une fonction complètement multiplicative vérifiant}
\begin{equation}
\forall p,\ |1+f(p)| \le 1
\end{equation}
\textit{et}
\begin{equation}
f(2) \not=-2\,.
\end{equation}

\textit{Alors les assertions $(7.7$) et $(7.8)$ suivantes sont équivalentes.}
\begin{equation}
\mathit{la\ fonction\ }f(n)/n\ \mathit{est\ }CMO
\end{equation}
\begin{equation}
\textit{pour\ tout\ r\'eel\ } \tau,\ \mathit{on\ a\ }\sum_p \frac{Re[(1+f(p))p^{-i\tau}]-1}{p} = -\infty\,.
\end{equation}

\vspace{2mm}

Ce théorème appelle plusieurs remarques.

\vspace{2mm}

\noindent\textbf{1)} Notons tout d'abord que l'on peut voir ce résultat de manière géométrique. Munissons l'ensemble des petites fonctions $CM$ de la norme
$$
\Big\|\frac{f(n)}{n}\Big\| = \sup_p |f(p)|\,.
$$
Alors le théorème 7 montre que toute petite fonction $CM$ qui est à une distance strictement inférieure à 1 de l'exemple d'Euler, est $CMO$. De plus, il décrit à quelles conditions une petite fonction $CM$ sur le cercle de rayon 1 est $CMO$.

\vspace{2mm}

\noindent \textbf{2)} Rappelons (voir la condition (2.2) du théorème~2) que la condition (7.6) est une condition nécessaire évidente pour que $f(n)/n$ soit $CMO$. En effet si $f(2)=-2$, on a $f(2^k)/2^k=(-1)^k$ et la suite $f(n)/n$ ne tend pas vers~0.

\vspace{2mm}

\noindent\textbf{3)} Ici on a
\begin{equation}
\forall p, \ Re \ f(p) \le 0\,.
\end{equation}

\noindent En choisissant $\tau=0$ dans (7.8), on retrouve donc la condition nécessaire
\begin{equation}
\sum_p \frac{Re\ f(p)}{p}= -\infty
\end{equation}
du théorème 6.

\vspace{2mm}

\noindent\textbf{4)} L'ensemble des fonctions multiplicatives est un groupe pour la convolution, alors que l'ensemble des fonctions complètement multiplicatives ne l'est pas. C'est la raison pour laquelle, pour prouver le théorème~7, nous serons amenés à évoluer dans le monde plus général des fonctions multiplicatives de somme nulle. En particulier, nous utiliserons le théorème~4.

Supposons à présent que l'on a toujours (7.9), mais que l'on remplace la condition~(7.5) par $\forall p$, $|f(p)|\le 1$. Alors l'assertion~(7.10) devient une condition nécessaire et suffisante pour que $f(n)/n$ soit~$CMO$.

\vspace{2mm}

\noindent \textbf{Théorème 8}.--- \textit{Soit $f$ une fonction complètement multiplicative vérifiant}
\begin{equation}
\forall p,\ |f(p)| \le 1
\end{equation}
\textit{et }
\begin{equation}
\forall p,\ Re\, f(p) \le 0\,.
\end{equation}

\textit{Alors les assertions $(7.13)$ et $(7.14)$ suivantes sont équivalentes}
\begin{equation}
\mathit{la\ fonction\ }f(n)/n\ \mathit{est\ } CMO		
\end{equation}
\begin{equation}
\sum_p \frac{Re\,f(p)}{p} = -\infty\,.
\end{equation}
\vspace{1mm}

Dans les grandes lignes, la preuve du théorème~8 est de même nature que celle du théorème~7. A la place du théorème de Hal\'asz, on utilise cette fois--ci le travail de majoration de la somme $\sum\limits_{n\le x}f(n)/n$ qu'ont effectué Montgomery et Vaughan en 2001 dans \cite{MoVa}, dans le cas où $f$ est une fonction complètement multiplicative à valeurs dans le disque unité.

Le travail de Hal\'asz, comme celui de Montgomery et Vaughan, repose sur une majoration uniforme sur des intervalles verticaux de longueur~1 situés au voisinage à droite de la droite $Re\,s=1$, de la fonction
$$
F(s)= \sum_n f(n)/n^s = \prod_p (1-f(p)/p^s)^{-1}\,.
$$
En majorant $F(s)$ sur les mêmes intervalles, mais cette fois--ci en moyenne, nous avons obtenu le résultat suivant.

\vspace{2mm}

\noindent \textbf{Théorème 9}.--- \textit{Soit $f$ une fonction complètement multiplicative vérifiant}
\begin{equation}
\forall p,\ |f(p)|\le 1
\end{equation}
\textit{et}
\begin{equation}
\limsup_{x\rg +\infty} \frac{\sum\limits_{p\le x}|1+f(p)|^2}{x/\log x} <1\,.
\end{equation}

\textit{Alors la fonction $f(n)/n$ est $CMO$.}

\vspace{4mm}
Nous avons donné aux théorèmes 7, 8 et 9, trois familles d'exemples de fonctions $CMO$, qui généralisent l'exemple historique d'Euler. Il serait intéressant d'unifier et généraliser ces familles d'exemples. Peut--on aller jusqu'à caractériser les petites fonctions $CMO$ parmi les petites fonctions~$CM$ ?

\section{Inversion de Möbius, théorèmes 10 et 11}
\setcounter{equation}{0}

On a donné au paragraphe précédent les grandes lignes de la preuve de l'implication $(7.8) \Rightarrow (7.7)$ du théorème~7. Pour établir la réciproque, nous avons besoin de relier le comportement asymptotique de $\sum\limits_{n\le x}f(n)/n$ à celui de $\sum\limits_{n\le x}g(n)$, sous une hypothèse de régularité plus générale que celle, (7.4), du théorème I+S. Nous donnons deux résultats : le premier généralise le théorème d'Ingham qui intervient dans le théorème I+S (\cite{Ing} ou chapitre II.7 de \cite{Ten}), le second celui de Ska\l ba (\cite{Ska} ou chapitre~II.7 de \cite{Ten}).

\vspace{2mm}

\noindent \textbf{Théorème 10}.--- \textit{Soient $f$ et $g$ deux fonctions arithmétiques vérifiant}
$$
g=f*1\,.
$$
\textit{On suppose qu'il existe un réel $\tau$ tel que}
\begin{equation}
\sum_{n\le x}g(n) = x^{1+i\tau}L(\log x) + o(x)\,,\quad (x\rg +\infty)
\end{equation}
\textit{où $L$ est une fonction bornée qui vérifie}
\begin{equation}
\sup_{x\le t\le x+1} |L(t) -L(x)| = o(1)\,,\qquad (x\rg +\infty)\,.
\end{equation}

\textit{Alors on a}
$$
\sum_{n\le x}\frac{f(n)}{n} = \frac{1}{x}\Big[\sum_{n\le x} f(n) + \Big(\frac{1}{i\tau\zeta(1+i\tau)}\Big)\sum_{n\le x}g(n)\Big] +o(1)\,,\quad (x\rg +\infty)\,.
$$
\textit{Si $\tau=0$, on convient que} $\frac{1}{i\tau\zeta(1+i\tau)}=1$.
\vspace{4mm}

Notons que le théorème d'Ingham correspond au cas où $\tau=0$ et $L$ est une fonction constante. Comme pour le résultat d'Ingham, le théorème~10 s'obtient en effectuant une sommation d'Abel suivie d'une inversion de Möbius qui utilise le théorème des nombres premiers.

\vspace{4mm}

\noindent \textbf{Théorème 11}.--- \textit{Soient $f$ et $g$ deux fonctions arithmétiques vérifiant}
$$
g=f*1\,.
$$
\textit{On suppose que $g$ est bornée et qu'il existe un réel $\tau$ tel~que}
\begin{equation}
\sum_{n\le x}g(n) = x^{1+i\tau}L(\log x) + o(x)\,,\quad (x\rg +\infty)
\end{equation}
\textit{où $L$ est une fonction qui vérifie}
\begin{equation}
\sup_{x\le t \le x+1} |L(t)-L(x)| = o(1)\,, \qquad (x\rg +\infty)
\end{equation}

\textit{Alors on a}
$$
\sum_{n\le x}f(n) = \frac{1}{\zeta(1+i\tau)} \sum_{n\le x}g(n) + o(x)\,,\quad (x\rg +\infty)\,.
$$
\textit{Si $\tau=0$, on convient que} $\frac{1}{\zeta(1)}=0$.
\vspace{4mm}

Remarquons que la fonction $L$ est nécessairement bornée. Notons par ailleurs que le théorème de Ska\l ba correspond au cas où $\tau=0$ et $L$ est une fonction constante. Là encore la preuve s'obtient en généralisant celle du résultat de Ska\l ba. Comme pour le théorème~10, on va effectuer une inversion de Möbius, à laquelle cette fois on ajoute le principe de l'hyperbole de Dirichlet, et on va utiliser le théorème des nombres premiers.

\section{Exemples d'Artin de fonctions $CMO$ ?}

Soient $L/K$ une extension galoisienne de corps de nombres,
\vspace{1mm}

\begin{tabular}{rl}
$G$ & $=\Gal(L/K)$ son groupe de Galois, \\ 
$\rho$ & une représentation de $G$ dans $GL_n(\C)$. 
\end{tabular} 

\vspace{1mm}
\noindent La fonction $L$ d'Artin $L(\rho,s)$ associée à cette représentation est alors, pour $Re\, s>1$, la somme d'une série de Dirichlet $L(\rho,s)=\sumn \frac{a_n}{n^s}$ où $(a_n)_{n\ge 1}$ est une suite multiplicative (pour tout ce paragraphe sur les fonctions $L$ d'Artin, voir \cite{Hei}). De plus, on sait que $L(\rho,s)$ admet un prolongement méromorphe sur $\C$. On a même la

\vspace{2mm}

\noindent\textbf{Conjecture d'Artin}

Si $\rho$ n'est pas la représentation triviale, alors $L(\rho,s)$ est une fonction entière.

\vspace{2mm}

Il est à noter en particulier que pour la plupart des séries $L$ d'Artin $L(\rho,s)=\sumn \frac{a_n}{n^s}$, on ne sait pas déterminer l'abscisse de convergence de cette série de Dirichlet.

Si $a_n$ est une fonction complètement multiplicative, on a $K=\Q$ et la fonction $L(\rho,s)$ est en réalité une fonction $L_\chi(s)$ associée à un caractère $\chi$ de Dirichlet. Considérons donc ici le cas où $a_n$ n'est pas complètement multiplicative.

Soit $\rho$ un nombre complexe tel que $\sumn a_n/n^\rho=0$. Comme $a_n$ n'est pas $CM$, $a_n/n^\rho$ n'est pas $CMO$. En revanche, on peut espérer qu'il existe une fonction $CMO$ \og au voisinage \fg{}  de~$a_n/n^\rho$.

De fait, s'il existe un $\rho=\beta +i \gamma$ avec $\beta>1/2$ tel~que
$$
\sumn \frac{a_n}{n^{\beta+i\gamma}} = 0\,,
$$
alors, en utilisant les propriétés de la suite $a_n$ et le théorème~4, on peut montrer que la fonction $CM$ $f(n)/n^\rho$ où $f$ est définie sur les nombres premiers par $f(p)=a_p$, est $CMO$. Il est à noter cependant que l'existence de ce type d'exemple de fonction $CMO$ est sujette à caution, dans la mesure où cela invaliderait l'hypothèse de Riemann généralisée aux fonctions $L$ d'Artin.

En réalité, si on suppose vérifiée l'hypothèse de Riemann généralisée aux fonctions $L$ d'Artin, nous ne  savons pas construire à partir de leurs zéros, d'autres fonctions $CMO$ que celles correspondant aux caractères $\chi$ de\break Dirichlet.

\section{Preuve du théorème 10}
\setcounter{equation}{0}

On note
$$
F(x) = \sum_{n\le x} f(n)
$$
et
$$
G(x) = \sum_{n\le x} g(n)\,.
$$
Par une sommation d'Abel et une inversion de Möbius, on obtient
\begin{align}
&\displaystyle\sum_{n\le x} \frac{f(n)}{n} - \frac{F(x)}{x}\nonumber\\
= & \int_1^x \frac{F(u)}{u^2}du\nonumber\\
= &\int_1^x \sum_{n\le u} \mu(n) G(u/n) \frac{du}{u^2}\nonumber\\
= &\sum_{n\le x} \mu(n)\int_n^x  G(u/n) \frac{du}{u^2}\nonumber\\
= &\sum_{n\le x} \frac{\mu(n)}{n}\int_1^{x/n}  G(v) \frac{dv}{v^2}\nonumber\\
=&\int_1^x \frac{G(v)}{v} \sum_{n\le x/v} \frac{\mu(n)}{n}\frac{dv}{v}\nonumber\\
=&P+E_1+E_2 \qquad \textrm{où}\nonumber
\end{align}
$$
P:= \int_1^x L(\log v) \sum_{n\le x/v} \frac{\mu(n)}{n}\frac{dv}{v^{1-i\tau}}\,,
$$
$$
E_1:= \int_1^{\sqrt{x}}\alpha(v) \sum_{n\le x/v} \frac{\mu(n)}{n} \frac{dv}{v}\ \mathrm{et}\ E_2:= \int_{\sqrt{x}}^x \alpha(v)\sum_{n\le x/v} \frac{\mu(n)}{n}\frac{dv}{v}
$$
avec
\begin{equation}
\alpha(v) := \frac{G(v)}{v} -v^{i\tau} L(\log v) =o(1)\,,\quad (v\rg +\infty)
\end{equation}

En utilisant le théorème des nombres premiers sous la forme 
\begin{equation}
\sum_{n\le v} \mu(n)/n \ll 1/\log^2 2v
\end{equation}
on obtient
$$
E_1 \ll \int_1^{\sqrt{x}} \frac{dv}{v\log^2 x/v} = \int_{\sqrt{x}}^x \frac{dv}{v\log^2 v} \ll \frac{1}{\log x}
$$
et
\begin{align}
E_2 &= \int_1^{\sqrt{x}} \alpha(x/v) \sum_{n\le v} \frac{\mu(n)}{n} \frac{dv}{v}\nonumber\\
&\ll \sup_{t\ge \sqrt{x}} |\alpha(t)| \cdot \int_1^{\sqrt{x}} \frac{dv}{v\log^2 2v} \newsmile \sup_{t\ge \sqrt{x}} |\alpha(t)|\,.\nonumber
\end{align}

Par ailleurs, on a
\begin{align}
P &= x^{i\tau} \int_1^x L(\log x/v) \sum_{n\le v} \frac{\mu(n)}{n} \frac{dv}{v^{1+i\tau}}\nonumber\\
&= x^{i\tau} L(\log x) Z_\tau(x) +O(E_3)\nonumber
\end{align}
\begin{align}
\mathrm{avec} \hspace{3cm} Z_\tau(x)&= \int_1^x \sum_{n\le v} \frac{\mu(n)}{n} \frac{dv}{v^{1+i\tau}}\nonumber\\
\mathrm{et} \hspace{4cm} E_3 &= \int_1^x \Big|L(\log x) -L(\log x/v)\Big| \ \Big|\sum_{n\le v} \frac{\mu(n)}{n}\Big| \frac{dv}{v}\,.\nonumber
\end{align}
Supposons que $\tau\not=0$. D'après l'exercice 178 de \cite{Ten}, on a
\begin{equation}
\sum_{n\le x} \frac{\mu(n)}{n^{1+i\tau}} = \frac{1}{\zeta(1+i\tau)} +O_\tau\Big(\frac{1}{\log x}\Big)\,.
\end{equation}
En utilisant également (10.2), on en déduit que
\begin{align}
Z_\tau(x) &=\sum_{n\le x} \frac{\mu(n)}{n} \int_n^x \frac{dv}{v^{1+i\tau}}\nonumber\\
&= \frac{1}{i\tau} \Big(\sum_{n\le x} \frac{\mu(n)}{n^{1+i\tau}}- x^{-i\tau} 
\sum_{n\le x} \frac{\mu(n)}{n}\Big)\nonumber\\
&= \frac{1}{i\tau\zeta(1+i\tau)} +O_\tau \Big(\frac{1}{\log x}\Big)\,.\nonumber
\end{align}
Par un calcul analogue, on obtient également
$$
Z_o(x)=1+O\Big(\frac{1}{\log x}\Big)\,.
$$

Notons 
\begin{equation}
\varphi(x) := \sup_{t\ge x} \sup_{1\le u\le e} |L(\log t)-L(\log t/u)|\,.
\end{equation}
En utilisant (10.2), on a uniformément pour $2\le y\le x$,
\begin{align}
E_3 &\ll \sup_{1\le u\le y} |L(\log x) -L(\log x/u)| \int_1^y \frac{dv}{v\log^2 2v} + \int_y^x \frac{dv}{v\log^2 2v}\nonumber\\
&\ll (\log y) \varphi(x/y) +1/\log y\,.\nonumber
\end{align}

Récapitulons. Avec ces différentes estimations, on a
$$
\begin{array}{c}
\displaystyle\sum_{n\le x} \frac{f(n)}{n} - \frac{F(x)}{x} - \frac{x^{-i\tau} L(\log x)}{i\tau\zeta(1+i\tau)}\\
\displaystyle\ll (\log y) \varphi(x/y) + \sup_{|t|\ge \sqrt{x}}|\alpha(t)| + \frac{1}{\log y} := \beta(x,y)\,.
\end{array}
$$

\noindent En choisissant $y=\min (\exp\Big\{\frac{1}{\sqrt{\varphi}(\sqrt{x})}\Big\},\ \sqrt{x})$ et en utilisant les hypothèses (8.1) et (8.2), on obtient que $\beta(x,y)=o(1)$ quand $x\rg +\infty$. En réutilisant  l'hypothèse (8.1) on conclut que
$$
\sum_{n\le x} \frac{f(n)}{n} = \frac{1}{x} \Big[ F(x) +\Big(\frac{1}{i\tau\zeta(1+i\tau)}\Big) G(x)\Big] +o(1)\,,\quad (x\rg +\infty)\,.
$$

\section{Preuve du théorème 11}

\noindent\textbf{Remarque préliminaire}

Au théorème 10, nous devions estimer $\int^x\frac{dF(u)}{u}$. Ici c'est $\int^x dF(u)$. Le poids 1 étant plus grand que le poids $1/u$, l'estimation ici est plus délicate. De fait on introduit ici l'hypothèse supplémentaire $g$ bornée. De plus, comme pour le théorème~10, on a recours au théorème des nombres premiers et à l'inversion de Möbius, mais on est amené ici à affiner en utilisant le principe de l'hyperbole de Dirichlet.

\vspace{2mm}

\noindent\textbf{Preuve.}

Nous utilisons les mêmes notations que pour le théorème~10. On a par le principe de l'hyperbole, pour tous $x>0$ et $y>0$,
$$
\begin{array}{c}
F(x) = \displaystyle \sum_{m\le y} \mu(m) G(x/m) + \sum_{n\le x/y} g(n)M(x/n) -G(x/y)M(y)\qquad \mathrm{avec}\\
\noalign{\vskip 2mm}
\displaystyle M(t):= \sum_{n\le t} \mu(n)\,.
\end{array}
$$
On utilise ici le théorème des nombres premiers sous la forme
$$
M(t) \ll t/\log^2 t\,, \qquad \qquad (t\ge 2)\,.
$$
Supposons dorénavant $x\ge 2y\ge 4$.

\noindent Comme $g$ est bornée, on a d'une part
$$
G(x/y)M(y) \ll x/\log^2 y
$$
et d'autre part
$$
\begin{array}{c}
\displaystyle \sum_{n\le x/y} g(n)M(x/n) \ll \sum_{n\le x/y} \frac{1}{n\log^2 (x/n)}\\
\noalign{\vskip2mm}
\displaystyle\ll x\int_{3/2}^{x/y} \frac{dt}{t\log^2x/t} = x \int_{y}^{2x/3} \frac{dt}{t\log^2t} \le \frac{x}{\log y}\,.
\end{array}
$$
On a enfin
$$
\sum_{m\le y} \mu(m) G(x/m) = x [x^{i\tau}(P-E_1)+E_2]
$$
avec
$$
\begin{array}{rl}
 P &:= L(\log x)\displaystyle \sum_{m\le y}\frac{\mu(m)}{m^{1+i\tau}}\,,\\
 \noalign{\vskip2mm}
 E_1 &:=\displaystyle \sum_{m\le y} \frac{\mu(m)}{m^{1+i\tau}}(L(\log x) -L(\log x/m))\\
\mathrm{et}\hspace{2cm}&\\
 E_2 &:=\displaystyle \sum_{m\le y} \frac{\mu(m)}{m} \alpha(x/m)\,.
\end{array}
$$
(On rappelle que la fonction $\alpha$ a été définie en (10.1)).

\noindent D'après (10.3), on a
$$
P=\frac{L(\log x)}{\zeta(1+i\tau)} \Big( 1+O_\tau \Big(\frac{1}{\log y}\Big)\Big)\,.
$$
Par ailleurs, on a avec $\varphi(t)$ définie en (10.4),
$$
E_1 \ll (\log^2 y) \varphi(x/y)\,.
$$
Enfin on a
$$ 
E_2 \ll (\log y) \sup_{t\ge x/y} |\alpha(t)|\,.
$$
Récapitulons. On a pour $x\ge 2y\ge 4$
$$
\frac{F(x)}{x^{1+i\tau}} = \frac{L(\log x)}{\zeta(1+i\tau)} + O\Big( (\log y) \sup_{t\ge x/y}|\alpha(t)| + (\log^2y) \varphi(x/y) +1/\log y\Big)\,.
$$
En utilisant (8.3) et (8.4) et en procédant comme à la conclusion de la preuve du théorème~10, on conclut~que
$$
F(x) = \frac{G(x)}{\zeta(1+i\tau)} +o(x)\,, \qquad (x\rg +\infty)\,.
$$

\section{Preuve du théorème 7}

\setcounter{equation}{0}

Notons $g=f*1$ et désignons par $\tilde{g}$ et $\tilde{f}$ les fonctions multiplicatives définies par 
$$
\tilde{g}(p^k) =g(p)\quad \mathrm{pour\ tout\ } p\ \mathrm{ premier\ et\ tout\ } k\ge 1
$$
et
$$
\tilde{g}= \tilde{f}*1 \,.
$$ 
D'après (7.5), on a $\|\tilde{g}\|_\infty \le 1$. On peut donc appliquer le théorème de Hal\'asz (théorème~III.4.5 de \cite{Ten}), ce qui nous amène à distinguer deux cas.

Si (7.8) est vérifiée, on a $\sum\limits_{n\le x}\tilde{g}(n)=o(x)$. En appliquant les théorèmes~10 et 11, on en déduit que
$$
\sumn \frac{\tilde{f}(n)}{n}=0\,.
$$

Supposons à présent que (7.8) n'est pas vérifiée. Alors pour le réel $\tau$ tel que $\sum\limits_p \frac{Re[\tilde{g}(p)p^{-i\tau}]-1}{p} > -\infty$, on~a
$$
\sum_{n\le x} \tilde{g}(n) = x^{1+i\tau} L(\log x) +o(x)
$$
où $L$ est une fonction à valeurs dans le cercle unité qui vérifie~(8.2).

\noindent En appliquant les théorèmes 10 et 11, on a donc
$$
\sum_{n\le x} \frac{\tilde{f}(n)}{n} = \Big(\frac{1}{i\tau\zeta(1+i\tau)} + \frac{1}{\zeta(1+i\tau)}\Big) x^{i\tau} L(\log x) +o(1)\,.
$$
Rappelons les conventions : si $\tau=0$, $1/\zeta(1+i\tau)=0$ et $1/(i\tau\zeta(1+i\tau))=1$. On~a~donc
$$
\frac{1}{i\tau\zeta(1+i\tau)} +\frac{1}{\zeta(1+i\tau)}\not=0 \quad \textrm{pour\ tout\ réel\ } \tau\,.
$$

On déduit de l'étude de ces deux cas que l'assertion~(7.8) est équivalente à l'identité
\begin{equation}
\sumn \frac{\tilde{f}(n)}{n}=0\,.
\end{equation}

Par ailleurs, en utilisant les hypothèses (7.5) et (7.6), on obtient, d'une part
$$
\forall p,\ \Big|\frac{f(p)}{p}\Big| <1\,,
$$
d'autre part
$$
\forall p,\ \sumk \frac{\tilde{f}(p^k)}{p^k} = 1+ \frac{f(p)}{p} \not= 0\,,
$$
et enfin
$$
\sum_p \Big|\big(1-\frac{f(p)}{p}\Big) \sumk \frac{\tilde{f}(p^k)}{p^k} -1 \Big| = \sum_p \Big| \frac{{f}(p)}{p} \Big|^2 < + \infty\,.
$$

\noindent En appliquant le théorème 4, on obtient donc que les assertions suivantes~(12.2) et (12.3) sont équivalentes
\begin{eqnarray}
\sumn \frac{\tilde{f}(n)}{n} = 0\ \\
\sumn \frac{f(n)}{n} =0\,.
\end{eqnarray}
On a donc prouvé que les assertions (7.8) et (12.3) sont équivalentes, ce qui conclut.

\section{Preuve du théorème 8}

\setcounter{equation}{0}

Compte tenu de l'hypothèse (7.12), l'implication de (7.13) vers (7.14) résulte du théorème~6. 

Attaquons nous maintenant à la réciproque (en réalité, l'hypothèse (7.14) n'interviendra qu'à la fin de notre argumentation). Notons comme d'habitude $F(s) =\sum\limits_n f(n)/n^s= \prod\limits_p (1-f(p)/p^s)^{-1}$, et 
$$
H(\sigma) := \Big(\sum_{k\in \Z} \max_{s=\sigma+it\atop |t-k|\le 1/2} \Big|\frac{F(s)}{s-1}\Big|^2\Big)^{1/2} := \Big(\sum_{k\in\Z} \max_{|t-k|\le 1/2} \Big|\frac{F(s)}{s-1}\Big|^2\Big)^{1/2}\,.
$$
Cette réciproque repose sur le résultat suivant, qui est implicite dans le travail \cite{MoVa} de Montgomery et Vaughan de~2001.

\vspace{2mm}

\noindent\textbf{Lemme fondamental $MV$}.---

 \textit{Soit $f$ une fonction $CM$ vérifiant $\sup\limits_p|f(p)|\le 1$ et $\lim\limits_{\sigma\searrow 1}(\sigma-1)H(\sigma)=0$. Alors la fonction $f(n)/n$ est~$CMO$. }

\vspace{4mm}

Cela résulte de la majoration
$$
\sum_{n\le x} \frac{f(n)}{n} \ll \frac{1}{\log x} \Big[H\Big(1+\frac{1}{\log x}\Big) + \int_{1+1/\log x}^2 \frac{H(\sigma)}{\sigma-1}d\sigma\Big]
$$
qui elle même constitue une légère variante du théorème~2 de \cite{MoVa}, et se démontre de la même manière.

Le travail effectif va consister dans un premier temps à majorer $F(s)$, pour pouvoir dans un deuxième temps majorer~$H(\sigma)$.

Pour tout réel positif $\tau$, notons 
$$
\tau^+ =\max(\tau,1) \qquad\mathrm{et}\qquad \tau^-=\min\,(\tau,1)\,.
$$

\noindent \textbf{Lemme 1}.--- \textit{Soit $f$ une fonction $CM$ vérifiant $\sup\limits_p |f(p)|\le 1$. Alors }
\begin{equation}
F(s) \ll F(\sigma)\Big(1+\frac{|t|}{\sigma-1}\Big)^{4/\pi}\,, \qquad (1<\sigma<2,t\in \R)\,.
\end{equation}
\textit{De plus, sous l'hypothèse supplémentaire $\sup\limits_p Re\, f(p)\le0$, il existe une cons\-tante $C_1>0$ telle que}
\begin{equation}
F(s) \ll \Big[1+ \frac{|t|^-}{\sigma-1}\Big]^{\frac{1}{2}+\frac{1}{\pi}} \ e^{C_1|t|^+} \quad (1<\sigma\le 2,t\in \R)\,.
\end{equation}

\vspace{2mm}

\noindent\textbf{Démonstration.} La majoration (13.1) est la première du lemma~2 de~\cite{MoVa}.

Pour prouver (13.2) on traite uniquement le cas où $t\ge 0$. Les calculs sont similaires pour $t\le 0$. 

On a pour~$\sigma>1$
$$
|F(s)| = \Big| \prod_p (1-f(p)p^{-s})^{-1}\Big| \newsmile \exp \Big\{\sum_p Re(p^{-it}f(p))/p^\sigma\Big\}\,,
$$
d'où, sous les hypothèses $\sup\limits_p |f(p)|\le 1$ et $\sup Re\ f(p)\le 0$,
\begin{equation}
F(s) \ll \exp \Big\{\sum_p h(t\log p)/p^\sigma\Big\},\quad (1<\sigma\le2,t\in \R)
\end{equation}
où $h(v) := \max\limits_{|z|\le 1\ \mathrm{et}\ Re\, z\le 0} Re(ze^{-iv})$.

En réalité, on observe facilement que $h$ est $2\pi$--périodique et vérifie
$$
h(v) = \left|
\begin{array}{cl}
|\sin v|	&\mathrm{si}\ |v|\le \pi /2\\
1 &\mathrm{si}\ \pi/2 \le v\le 3\pi/2
\end{array}
\right.
$$
Par une sommation d'Abel, on a
$$
\begin{array}{rl}
\displaystyle\sum_p \frac{h(t\log p)}{p^\sigma} &\displaystyle=\int_e^{+\infty} \frac{h(t\log u)}{u^\sigma\log u} du \ + O(t^+)\\
&\displaystyle=\int_t^{+\infty} h(v) e^{-v(\sigma-1)/t}\ \frac{dv}{v}\ + O(t^+)\,.
\end{array}
$$
D'où
\begin{equation}
\sum_p  \frac{h(t\log p)}{p^\sigma} = Q(t) +O(t^+)
\end{equation}
avec
$$
Q(t) := \int_{t^+}^{+\infty} h(v) e^{-v(\sigma-1)/t}\ \frac{dv}{v}\,.
$$

Si $t\le \sigma-1$, alors $Q(t)=O(1)$, ce qui avec (13.3) établit la majoration~(13.2).

Supposons dorénavant $t>\sigma-1$. On a alors
$$
\begin{array}{rl}
Q(t) &=\displaystyle \int_{t^+}^{\frac{t}{\sigma-1}} h(v) \frac{dv}{v} +O(1)\\
\noalign{\vskip2mm}
&= \displaystyle\sum_{t^+/2\pi\le k<t/2\pi(\sigma-1)} \int_{2\pi k}^{2\pi(k+1)} h(v) \Big(\frac{1}{2\pi k} +O\Big(\frac{1}{k^2}\Big)\Big) dv+O(1)\\
\noalign{\vskip2mm}
&\displaystyle=\frac{1}{2\pi} \Big(\int_0^{2\pi}h(v)dv\Big)\log \Big(\frac{t}{t^+(\sigma-1)}\Big) +O(1)\\
\noalign{\vskip2mm}
&=\displaystyle\Big(\frac{1}{2}+\frac{1}{\pi}\Big) \log \Big(\frac{t^-}{\sigma-1}\Big) +O(1)\,.
\end{array}
$$
Donc là encore, la majoration (13.2) découle grâce à la combinaison de (13.3) et~(13.4).

Nous aurons également besoin de la majoration
\begin{equation}
|F(s)| \le \zeta(\sigma) \newsmile \frac{1}{\sigma-1}, \ (1<\sigma\le 2)\,.
\end{equation}

En utilisant successivement les majorations (13.5) et (13.2) on obtient, pour la constante $C_1$ du lemme~1,
$$
\sum_{|k|\ge\frac{1}{6C_1}\log (\frac{2}{\sigma-1})} \max_{|t-k|\le \frac{1}{2}} \Big|\frac{F(s)}{s-1}\Big|^2 \ll \frac{1}{(\sigma-1)^2\log(\frac{2}{\sigma-1})}
$$
et 
$$
\sum_{1\le|k|<\frac{1}{6C_1}\log (\frac{2}{\sigma-1})} \max_{|t-k|\le \frac{1}{2}} \Big|\frac{F(s)}{s-1}\Big|^2 \ll \frac{1}{(\sigma-1)^{1+\frac{2}{\pi}+\frac{1}{3}}}\,.
$$
Donc 
\begin{equation}
\Big(\sum_{k\in\Z\ba\{0\}}\max_{|t-k|\le\frac{1}{2}} \Big|\frac{F(s)}{s-1}\Big|^2 \Big)^{1/2} \ll \frac{1}{(\sigma-1)\sqrt{\log\frac{2}{\sigma-1}}}\,.
\end{equation}

Supposons à présent $0\le t\le 1/2$ (les calculs sont identiques pour $-1/2\le t\le 0$). On note $a:=\frac{2\pi}{\pi+6}$.

Si $\frac{\sigma-1}{t} \le (F(\sigma))^a$, en utilisant la majoration (13.2), on~obtient
$$
\frac{F(s)}{s-1} \newsmile \frac{F(s)}{t} \ll \frac{1}{\sigma-1}\Big(\frac{\sigma-1}{t}\Big)^{\frac{1}{2}-\frac{1}{\pi}} \le \frac{1}{\sigma-1} (F(\sigma))^{a(\frac{1}{2}-\frac{1}{\pi})}\,.
$$
Si au contraire $(F(\sigma))^a\le \frac{\sigma-1}{t}$, en utilisant la majoration (13.1) on~obtient
$$
\frac{F(s)}{s-1} \ll \frac{F(s)}{\sigma-1} \ll \frac{F(\sigma)}{\sigma-1}\Big(1+\frac{|t|}{\sigma-1}\Big)^{\frac{4}{\pi}} \ll \frac{1}{\sigma-1} (F(\sigma))^{1-\frac{4a}{\pi}}\,.
$$
Or $a\Big(\frac{1}{2}-\frac{1}{\pi}\Big)=1-\frac{4a}{\pi}= 1-\frac{8}{\pi+6}$. En combinant ces deux cas, on obtient donc
\begin{equation}
 \max_{|t|\le 1/2} \Big|\frac{F(s)}{s-1}\Big| \ll \frac{F(\sigma)^{1-\frac{8}{\pi+6}}}{\sigma-1}\,.
 \end{equation}

Nous sommes maintenant en mesure de conclure. Soit $f$ une fonction $CM$ vérifiant (7.11) et (7.12). Alors
$$
F(\sigma) \newsmile \exp \Big\{ \sum_p \frac{Re\,f(p)}{p^\sigma}\Big\}\,, \quad (1<\sigma\le2)\,,
$$
et par convergence monotone,
$$
\sum_p \frac{Re\,f(p)}{p}	= \lim_{\sigma\searrow 1} \sum_p \frac{Re\, f(p)}{p^\sigma}\,.
$$
D'après l'hypothèse (7.14), on a donc
$$
\lim_{\sigma\searrow 1} F(\sigma)=0\,.
$$
Avec (13.6) et (13.7), on en déduit que
$$
\lim_{\sigma\searrow 1} (\sigma-1)H(\sigma)=0\,.
$$
Le lemme fondamental $MV$ permet donc de conclure que $f(n)/n$ est~$CMO$.

\section{Preuve du théorème 9}

\setcounter{equation}{0}

A toute série de Dirichlet $\sum\limits_{n\ge1} a_n n^{it}$, on associe la fonction 
$$
S(x)=|\{n\le x : a_n\not=0\}|\,.
$$
Nous aurons besoin du résultat auxiliaire suivant

\vspace{2mm}

\noindent \textbf{Lemme 2}.--- \textit{Il existe une constante $C_2>0$ telle que pour toute série de Dirichlet absolument convergente $\sum\limits_n a_n n^{it}$, on~a}
$$
\sup_{\tau\in \R}
 \int_{\tau-1/2}^{\tau+1/2} \Big| \sum_n a_n n^{it}\Big|^2 dt \le C_2 \sum_n |a_n|^2 \Big(S(9n)-S\Big(\frac{n}{9}\Big)\Big)\,.
 $$
 
La preuve s'obtient sans difficulté en généralisant celle de la formule~(1.4.12) du theorem~1.4.6 de~\cite{QuQu}.

\vspace{2mm}

\noindent\textbf{Remarque préliminaire.}

Dans le cas de l'exemple d'Euler $\lambda(n)/n$, on a $F(s)=\zeta(2s)/\zeta(s)$ qui est holomorphe au voisinage de la droite $Re\, s=1$. Cela n'est plus vrai dans le cas général. L'un de nos premiers objectifs est de montrer que les hypothèses du théorème~9 entraînent l'existence presque partout en $t$ de la limite $\lim\limits_{\sigma\rg 1^+}F(\sigma+it)$. Pour cela on travaille~avec
$$
G(s) := F(s)\zeta(s)
$$
et on va majorer en moyenne la fonction $G(s)$.

\vspace{2mm}

\noindent \textbf{Preuve du théorème 9}.

D'après l'hypothèse (7.15), les fonctions $F(s)$ et $G(s)$ sont holomorphes dans le demi--plan~$Re\, s>1$.

Soit $I$ un intervalle de $\R$ de mesure 1. On a pour tous $1<\sigma_1\le \sigma_2\le 2$
\begin{equation}
\int_I |G(\sigma_2+it) -G(\sigma_1+it)| dt \le \int_{\sigma_1}^{\sigma_2} \int_I |G'(\sigma+it)| dt\, d\sigma
\end{equation}
et par Cauchy--Schwarz
\begin{equation}
\Big(\int_I |G'(\sigma+it)|dt\Big)^2 \le J(\sigma) K(\sigma)
\end{equation}
avec
$$
J(\sigma) = \int_I \Big|\frac{G'}{G} (\sigma+it)\Big|^2 dt
$$
et
$$
K(\sigma) = \int_I |G(\sigma+it)|^2 dt\,.
$$

Notons
$$
G(s) = \sumn \frac{g(n)}{n^s}\,.
$$
La fonction $g(n)$ est multiplicative et vérifie
\begin{equation}
g(p^k) = \sum_{j=0}^k f^j(p)\,.
\end{equation}
Avec l'hypothèse $\|f\|_\infty\le 1$, on a donc pour $Re\, s>1$
$$
\frac{G'}{G}(s) = -\sum_p \frac{(\log p)g(p)}{p^s} +O(1)\,.
$$
D'où, en utilisant le lemme 2 et le théorème des nombres premiers,
$$
J(\sigma) \ll 1+ \sum_p \frac{\log p}{p^{2\sigma-1}} \le 1+ \sum_p \frac{\log p}{p^\sigma} \newsmile \frac{\zeta'}{\zeta}(\sigma) \newsmile \frac{1}{\sigma-1}\,.
$$
Par ailleurs en réutilisant le lemme 2, on a
$$
\begin{array}{rl}
K(\sigma) 	&\ll\displaystyle\sum_n \frac{|g(n)|^2}{n^{2\sigma-1}} \le \sum_n \frac{|g(n)|^2}{n^\sigma} = \prod_p \sum_{k\ge 0} \frac{|g(p^k)|^2}{p^{k\sigma}}\\
&=\displaystyle \exp \Big\{\sum_p \log\Big( 1+ \frac{|g(p)|^2}{p^\sigma} +O\Big(\frac{1}{p^2}\Big)\Big)\Big\}
\end{array}
$$
d'après (14.3) et l'hypothèse $\|f\|_\infty \le 1$. On a donc
$$
K(\sigma) \ll \exp \Big\{\sum_p \frac{|g(p)|^2}{p^\sigma}\Big\}\,.
$$
D'où avec (14.1) et (14.2),
\begin{equation}
\int_I |G(\sigma_2+it) -G(\sigma_1+it)|dt \ll \int_{\sigma_1}^{\sigma_2} \exp\Big\{\frac{1}{2} \sum_p \frac{|g(p)|^2}{p^\sigma}\Big\} \frac{d\sigma}{\sqrt{\sigma-1}}\,.
\end{equation}

Notons $A(t) = \sum\limits_{p\le t} |g(p)|^2$. On a
$$
\sum_p \frac{|g(p)|^2}{p^\sigma} =\sigma \int_2^{+\infty} \frac{A(t)}{t^{\sigma+1}} dt
$$
Or $g(p)=1+f(p)$. Donc l'hypothèse (7.16) du théorème~9 entraîne qu'il existe deux constantes $t_o>1$ et $c>0$ telles que
$$
\int_{t_o}^{+\infty} \frac{A(t)}{t^{\sigma+1}} dt \le (1-c) \int_{t_o}^{+\infty}
\frac{dt}{t^\sigma\log t}= (1-c)\log \frac{1}{\sigma-1} +O(1)\,.
$$
Il découle donc de la majoration (14.4) que
\begin{equation}
\sup_{\tau\in \R} \int_{\tau-1/2}^{\tau+1/2} |G(\sigma_2+it)-G(\sigma_1+it)| dt \le\! \int_{\sigma_1}^{\sigma_2} M(\sigma)d\sigma\,,\ (1\!<\!\sigma_1\!\le\! \sigma_2\!\le\! 2)
\end{equation}
où $M$ est une certaine fonction positive et intégrable sur $]1,2]$. Donc quand $\sigma\rg 1^+$, $G(\sigma+it)$ converge dans $L_{loc}^1(\R)$ vers sa limite ponctuelle qui est définie presque partout en $t$ et que l'on notera $G(1+it)$. La fonction\break $F(1+it):= \lim\limits_{\sigma\searrow 1}F(\sigma+it)$ est donc également définie pour presque tout réel~$t$.

Soient $x$ et $a$ deux réels fixés pour l'instant et vérifiant $x\ge 2\ge 2a>0$. Soit $\sigma\in ]1,2]$. On note
$$
\Phi_\sigma(t) = F(\sigma+it) \frac{\sin xt}{t} \Big(\frac{\sin at}{at}\Big)^2 = \frac{\sin xt}{t\zeta(\sigma +it)}G(\sigma+it)\Big(\frac{\sin at}{at}\Big)^2\,.
$$
D'après les propriétés (5.3) et (5.4) de la fonction $\zeta(s)$, on~a
$$
\lim_{\sigma\searrow 1\atop L^\infty(\R)} \frac{\sin xt}{t\zeta(\sigma+it)} = \frac{\sin xt}{t\zeta(1+it)}\,.
$$
Par ailleurs, on a prouvé que
$$
\lim_{\sigma\searrow 1\atop L_{loc}^1(\R)} G(\sigma+it) = G(1+it)\,.
$$
Compte tenu de la majoration (14.5), on en déduit que
$$
\lim_{\sigma\searrow 1\atop L^1(\R)} G(\sigma+it) \Big(\frac{\sin at}{at}\Big)^2 = G(1+it) \Big(\frac{\sin at}{at}\Big)^2\,.
$$
Donc finalement, on a
$$
\lim_{\sigma\searrow 1\atop L^1(\R)}\Phi_\sigma=\Phi_1\,.
$$
En particulier, on a
\begin{equation}
\lim_{\sigma\searrow 1^+}\int_\R F(\sigma+it) \frac{\sin xt}{t}\Big(\frac{\sin at}{at}\Big)^2 dt = \int_\R F(1+it)\frac{\sin xt}{t}\Big(\frac{\sin at}{at}\Big)^2dt\,.
\end{equation}

On peut à présent suivre la démarche du paragraphe 5. On a d'après (5.2), pour tout~$\sigma>1$,
$$
\begin{array}{l}
\displaystyle\frac{1}{\pi a} \int_\R \frac{F(\sigma +it)}{t} \Big(\frac{\sin at }{at}\Big)^2 \sin xt\, dt = \sum_n \frac{f(n)}{n^\sigma} w_{x,a}(\log n)\\
\displaystyle = \sum_{\log n\le x} \frac{f(n)}{n^\sigma} +O \Big(\sum_{x-2a\le \log n \le x +2a} \frac{1}{n}\Big) = \sum_{\log n\le x} \frac{f(n)}{n^\sigma} +O(a+e^{-x})\,.
\end{array}
$$
En passant à la limite quand $\sigma \rg 1^+$ et en utilisant (14.6), on~obtient
\begin{equation}
\sum_{\log n\le x} \frac{f(n)}{n} = \frac{1}{\pi a} \int_\R \frac{F(1+it)}{t} \Big(\frac{\sin at}{at}\Big)^2 \sin xt\, dt +O(a+e^{-x})\,.
\end{equation}
Or
$$
\begin{array}{l}
\displaystyle\int_\R \Big|\frac{F(1+it)}{t}\Big(\frac{\sin at}{at}\Big)^2\Big| dt \ll \int_\R \frac{|G(1+it)|}{1+t^2} dt\\\noalign{\vskip4mm}
\displaystyle\ll \sum_{k\in \Z} \frac{1}{k^2} \int_{k-1/2}^{k+1/2} |G(1+it)|dt < +\infty \quad \textrm{d'après}\ (14.5)\,.
\end{array}
$$
Donc $\frac{F(1+it)}{t}\Big(\frac{\sin at}{at}\Big)^2 \in L^1(\R)$.
En faisant tendre $x$ vers $+\infty$ dans (14.7) et en appliquant le lemme de Riemann--Lebesgue, puis en faisant tendre $a$ vers $0$, on conclut que $f(n)/n$ est~$CMO$.

\vskip4mm

\noindent\textbf{Remarque.}
\vspace{2mm}

Nous avons exprimé l'hypothèse (7.16) du théorème~9 sous cette forme pour avoir une condition suffisante sur la distribution des $f(p)$ qui s'exprime simplement. Mais on vérifie facilement que l'on a en réalité démontré le résultat légèrement plus fort suivant.

\vspace{4mm}

\noindent \textbf{Théorème $9'$}.--- \textit{Soit $f(n)$ une fonction complètement multiplicative vérifiant}
$$
\sup_p|f(p)| \le 1
$$
et
$$
\int_1^2 \exp \Big\{ \frac{1}{2} \sum_p \frac{|1+f(p)|^2}{p^\sigma}\Big\} \frac{d\sigma}{\sqrt{\sigma-1}} <+\infty\,.
$$
\textit{Alors $f(n)/n$ est $CMO$.}
\vspace{2mm}
\vspace*{1cm}

\noindent\textbf{Remerciements.}
\vspace{2mm}

C'est Michel Balazard qui est à l'origine de ce présent travail. Il nous a fait prendre conscience du théorème~1 en mettant en avant les exemples d'Euler et de Dirichlet--Riemann. Cela nous a donné envie d'en savoir plus. Nous le remercions pour tout~cela.



\vskip4mm

\begin{tabular}{ll}

Jean--Pierre Kahane & Eric Saias \\

Laboratoire de mathématique &Laboratoire de Probabilités et\\

Bâtiment 425& Modèles Aléatoires\\

Université Paris--Sud &Université Pierre et Marie Curie\\

91405 Orsay Cedex (France)  &4, place Jussieu\\

& 75252 Paris Cedex 05 (France)\\

\vspace{2mm}

\textsf{jean-pierre.kahane@u-psud.fr} &\textsf{eric.saias@upmc.fr}
\end{tabular}

\end{document}